\documentstyle[12pt]{article}
\renewcommand{\baselinestretch}{1.3}

\newtheorem {th}{Theorem}
\newtheorem {lem}[th]{Lemma}
\newtheorem {pr}[th]{Proposition}

\newtheorem{defn}{Definition}

\newtheorem{conj}[th]{Conjecture}

\def\Cox{\hfill \Box}
\def\mk{{\kern-.6em}}

\def\deq{\, {\stackrel {def} {=}}}

\def\dd{\delta}

\def\iid{\mbox{independent, identically distributed}}

\def\ee{\epsilon}
\def\ft{{\tilde{f}}}
\def\E{{\bf{E}}}
\def\P{{\bf{P}}}
\def\N{\hbox{I\kern-.2em\hbox{N}}}
\def\R{\hbox{I\kern-.2em\hbox{R}}}

\def\|{\, | \, }

\def\0{\hat{0}}
\def\1{\hat{1}}

\begin{document}
 
\begin{titlepage}
\begin{center}
{\large \bf DOMINATION BETWEEN TREES AND APPLICATION TO AN EXPLOSION
PROBLEM} \\
Robin Pemantle \footnote{Oregon State University and the University 
of Wisconsin-Madison, Dept. of Math, Van Vleck Hall, 480 Lincoln
Drive, Madison, WI 53706}$^,$
\footnote{Research supported in part by National Science Foundation 
Grant \# DMS 9103738}
and Yuval Peres \footnote{Dept. of Math, Yale University, New Haven, 
CT 06520 . }
\end{center}

\vfill

{\bf ABSTRACT:} \break
We define a notion of stochastic domination between trees, where one
tree dominates another if when the vertices of each are labeled
with $\iid$ random variables, one tree is always more likely to contain
a path with a specified property.  Sufficient conditions for this
kind of domination are (1) more symmetry and (2) earlier branching.
We apply these conditions to the problem of determining how fast
a tree must grow before first-passage percolation on the tree
exhibits an explosion, that is to say, infinitely many vertices
are reached in finite time.  For a tree in which each vertex at distance 
$n-1$ from the root has $f(n)$ offspring, $f$ nondecreasing,
an explosion occurs with exponentially distributed passage times
if and only if $\sum f(n)^{-1} < \infty$. 
\vfill

\noindent{Keywords:} first-passage percolation, explosion, tree,
domination

\noindent{Subject classification: } Primary: 60K35 ; Secondary: 
60F15 , 60E15.

\end{titlepage}

\section{Introduction}

Let $\Gamma$ be any locally finite tree with some vertex $\rho$ chosen 
as the root and total height $N \leq \infty$.  Label the vertices of 
$\Gamma$ with $\iid$ real random variables $\{ X(\sigma) : \sigma \
\in \Gamma \}$, and let $B \subseteq \R^N$ be some Borel set.  Let 
$P(B;\Gamma)$ denote the probability that $(X(\sigma_1) , X(\sigma_2)
,\ldots) \in B$ for some for some non-self-intersecting
path $\rho , v_1 , \ldots , v_N$ in $\Gamma$.  This probability
arises in many contexts.  The problem as such is studied in
\cite{Ly} and \cite{Ev}.  For first-passage percolation, the probability
of reaching level $N$ by time $T$ is $P(B;\Gamma_{|N})$, where $B$ is
the set $\sum_{i=1}^N v_i \leq T$ and $\Gamma_{|N}$ is the first $N$ levels
of $\Gamma$.  The same quantity arises when studying diffusion-limited
aggregation on trees via the exponential representation \cite{AS,BPP}.
A random walk in a random environment on a tree will be transient
when some path is itself transient, which probability can be reduced 
via the electrical representation to $P(B; \Gamma)$ for $\iid$
resistances, the set $B$ being again the summable sequences;
see \cite{Pe} and \cite{PP}.  The
study of tree-indexed Markov chains \cite{BP2} can be reduced to
computations of $P(B;\Gamma)$ by representing the Markov chain as
a function of $\iid$ uniform $[0,1]$ random variables.  

Let $|v|$ denote the distance from $\rho$ to $v$ and write $w \leq v$
if $w$ is on the path from $\rho$ to $v$.  Let $\Gamma_n
= \{ v : |v| = n \}$ denote the $n^{th}$ level of $\Gamma$
and $\Gamma_{|n}$ denote the first $n$ levels of $\Gamma$.
Say that $\Gamma$ is {\em spherically symmetric} with growth
function $f$ if every vertex $v \in \Gamma_n$ has $f(|v|+1)$
children (a child of $v$ is a neighbor $w$ with $v \leq w$).  
The notion that we examine in this paper of stochastic domination 
between trees is as follows.

\begin{defn}
Let $\Gamma$ and $\Gamma'$ be two finite or infinite trees with roots
$\rho$ and $\rho'$ respectively.  We say that $\Gamma$ {\em dominates}
$\Gamma'$ if whenever
$$\{ X(v) : v \in \Gamma \} \cup \{ X' (w) : w \in \Gamma' \} $$
is a collection of i.i.d. random variables, $n \geq 1$ is an integer,
and $B \subseteq \R^n$ is a Borel set, then
$$P (B ; \Gamma_{|n}) \geq P(B ; \Gamma_{|n}') .$$
\end{defn}

Spherically symmetric trees are easier to compute with, which
is the principal motivation for developing inequalities that compare
$P(B; \Gamma)$ and $P(B;\Gamma')$ when one of $\Gamma$ or $\Gamma'$
is spherically symmetric.  Our main comparison result, proved in
Section~2, is the following.

\begin{th} \label{comparison}
Let $\Gamma$ and $\Gamma'$ be finite or infinite trees with $\Gamma$
spherically symmetric.  Then $\Gamma$ dominates $\Gamma'$ if and only
if for every $n \geq 1$ the $n^{th}$ generation 
$$\Gamma_n \deq \{ v \in \Gamma : |v| = n \}$$
of $\Gamma$ is at least as big as the $n^{th}$ generation of $\Gamma'$.
\end{th}

{\bf Remarks:} 

\noindent{1. } We understand the domination partial order completely
only for spherically symmetric trees and for trees of height two.  In
the latter case, it reduces to the classical notion of Hardy majorization
(Proposition~\ref{height 2}). 

\noindent{2. } Comparison results for $P(B;\Gamma)$ as $B$ varies
may be found in \cite{PP}.

A consequence of Theorem~\ref{comparison} is that among all
trees $\Gamma$ of height $n$ having $|\Gamma_n| = k$, the
tree $T(n,k)$ consisting of $k$ disjoint paths of length
$n$ joined at the root is maximal in the domination order.
If the common law of the $X(\sigma)$ is $\mu$ and $B \subseteq \R^n$,
then $1 - P(B;T(n,k)) = (1 - \mu^n (B))^k$, thus for any $\Gamma$
of height $n$,
\begin{equation} \label{stringy}
1 - \P (B ; \Gamma) \geq \left ( 1 - \mu^n (B) \right )^k .
\end{equation}
The definition of $P(B;\Gamma)$ extends naturally to any
{\em graded graph}, this being a finite graph whose vertices
are partitioned into levels $1 , \ldots , n$ with oriented edges 
allowed only from levels $i$ to $i+1$, $i = 1 , \ldots , n-1$.
A natural conjecture is then
\begin{conj} \label{graded graphs}
If $G$ is a graded graph of height $n$, let $K(G)$ be
the number of oriented paths that pass through every level
of $G$.  Let $X(\sigma)$ be $\iid$ random variables with 
common law $\mu$ and $B \subseteq \R^n$.  Then
$$1 - \P (B ; G) \geq \left ( 1 - \mu^n (B) \right )^{K(G)} .$$
\end{conj}
If $B$ is an upwardly closed set then both the conjecture
and the inequality~(\ref{stringy}) follows easily from
the FKG inequality.  In the case where $n = 2$ ($G$ is
bipartite), Conjecture~\ref{graded graphs} is due to
Sidorenko \cite[Conjecture 5.2]{Si3} in the form of an 
analytic inequality; Sidorenko has proved several related
analytic inequalities on graphs in~\cite{Si1,Si2} and~\cite{Si3}
including the conjecture itself for the special cases where
$G$ is bipartite, acyclic, a single cycle or sufficiently small.

The remainder of the paper is devoted to an in-depth application of 
this theorem to first-passage percolation, which we now describe.  
Let $\Gamma$ be a spherically symmetric infinite tree with
growth function $f$, and let $\{ X(v) \}$ be a collection 
of independent exponential random variables of mean one.
(Our results hold for a much more general class distributions 
described later.) Think of $X(v)$ as a transit time
across the edge connecting $v$ to its parent.  Define
\begin{eqnarray*}
S(v) = \sum_{\rho < w \leq v} X(w) & \mbox{ and } & M_n = \min_{|v| = n} S(v) .
\end{eqnarray*}
In the context of first-passage percolation, $S(v)$ is the passage time 
from the root to the vertex $v$ and $M_n$ is the first passage time to
the $n^{th}$ generation of $\Gamma$.  Say that 
an explosion occurs if the increasing sequence $M_n$ is bounded.  

An explosion is a tail event, so for a given tree the probability
of an explosion is zero or one.  It is natural to try to determine when
the probability is zero and when it is one.  This problem was brought
to our attention by Enrique Andjel (personal communication) in 
reference to uniqueness proofs for particle systems via graphical
representations.  For arbitrary trees,
such questions are difficult to settle without leaving open some
``critical'' case (see for example \cite{LP}).  As observed by Andjel,
it is elementary, for the case of spherically symmetric trees,
to obtain conditions for explosion or non-explosion
which are almost sharp.  Indeed, if 
$$\liminf_{n \rightarrow \infty} (1/n!) \prod_{i=1}^n f(i) < \infty ,$$
then a simple Borel-Cantelli argument shows there is almost surely no
explosion before time one, from which it follows easily that with
probability one, no explosion occurs at all; on the other
hand, the condition $\sum f(n)^{-1} < \infty$ implies
that the greedy algorithm finds an explosion.  Thus in particular
there is no explosion when $f(n) = n$ but there is one when $f(n) =
n \ln (n)^{1+\ee}$.  We show that the summability condition is
sharp, the statement being a little more complicated in the
case where $f$ is not an increasing function.  

\begin{th} \label{increasing}
Let $\Gamma$ be a spherically symmetric tree with growth
function $f$ that is nondecreasing.  The probability
of an explosion is zero or one according to whether
$\sum_{n=1}^\infty f(n)^{-1}$ is infinite or finite.  
Furthermore, if the sum is infinite and $f(n)$ is unbounded as well,
then $M_N / \sum_{n=1}^N f(n)^{-1}$ converges almost surely to $e^{-1}$.  
\end{th}

In most integral tests in probability theory, some regularity condition
is imposed.  An unusual feature of the problem considered here is that
it permits a criterion (Theorem~\ref{general})
valid for arbitrary growth functions $f$.  
The condition in Theorem~\ref{increasing} always suffices for
explosion but simple examples ($f(2n-1) = 1 , f(2n) = 2^n$) show it
is not necessary.  For any $f : \N \rightarrow \N$, define a function
$\ft : \N \rightarrow \R^+$ recursively by
\begin{equation} \label{rec def}
\ft (n+1) = \sup \left \{ a : \;  a^m \prod_{i=1}^n 
   \ft (i) \leq \prod_{i=1}^{n+m} f(i) \mbox{ for all } m \geq 1 \right \} .
\end{equation}
In particular, 
$$\ft (1) = \inf_{m \geq 1} \left ( \prod_{i=1}^m f(i) \right )^{1/m} .$$
It is easy to see that the function $\ft$ is always nondecreasing and
coincides with $f$ when $f$ is nondecreasing.  

\begin{th} \label{general}
Let $\Gamma$ be a spherically symmetric tree with growth
function $f$, labeled as before by independent exponentials $X(v)$
of mean 1.  Then the probability of an explosion is zero or
one according to whether $\sum_{n=1}^\infty \ft (n)^{-1}$ is
infinite or finite.  Furthermore, if the sum is infinite and $\ft$
is unbounded then 
$$\limsup_{N \rightarrow \infty} M_N / \sum_{n=1}^N \ft (n)^{-1}
= e^{-1}.$$
\end{th} 

\noindent{{\bf Remarks}}

\noindent{1. } It is usually easy to calculate $\ft$ from $f$.  Informally,
if you graph $\sum_{j=1}^n \ln (f(j))$ against $n$, and consider the
convex hull of the region above this graph, its boundary is the graph of
$\sum_{j=1}^n \ln \ft (j)$.  For instance if the two sequences 
$\{ f(2n-1) \}$ and $\{ f(2n) \}$ are nondecreasing with
$f(2n-1) < f(2n)$ for all $n$, then 
$$\ft(2n-1) = \ft (2n) = \left ( f(2n-1) f(2n) \right )^{1/2} .$$

\noindent{2. } Theorem~\ref{general} is proved by comparing the tree
$\Gamma$ to a tree with nondecreasing growth function, which is
where the application of Theorem~\ref{comparison} is needed.

These theorems are proved in Section~3.  Versions where the
variables $\{ X(v)\}$ are not exponential are given in Section~4.
Results completely analogous to
Theorems~\ref{increasing} and~\ref{general} hold for distributions
$G$ satisfying a power law near zero, i.e. $\lim_{ t \downarrow 0}
G(t) / t^{\alpha}$ is finite and positive for some $\alpha > 0$.
Obtaining sharp criteria for explosion to occur that are valid for
arbitrary transit-time distributions $G$ seems more delicate, although
we cannot find a $G$ for which the natural criterion~(\ref{criterion})
in the last section fails.

\section{Domination between trees}

\noindent{{\bf Proof of Theorem~{\protect \ref{comparison}}}}: We start 
by establishing
the theorem under the additional assumption that $\Gamma'$ is
spherically symmetric.  This is the only case that is used
in the proof of Theorem~\ref{general}.

Let $\mu$ be any probability measure on $\R$ and let $D$ be a Borel
set in $\R^n$.  If $b_1 , \ldots , b_n$ are the cardinalities 
of the generations of a spherically symmetric tree (i.e. the
growth function is $f(i) = b_i / b_{i-1}$), let the vertices 
index i.i.d. random variables $X(v)$ with common law $\mu$ and
let $\Psi (b_1 , \ldots , b_n ; D)$ denote the probability
that all paths $\rho , v_1 , \ldots , v_n$ in the tree satisfy
$(X(v_1) , \ldots , X(v_n)) \in D$.  Passing to complements in
the definition of domination, we must show that
$$\Psi (b_1 , \ldots , b_n ; D) \leq \Psi (b_1' , \ldots , b_n' ;D)$$
whenever both are defined and $b_i \geq b_i'$ for all $i$.  

The key to doing this is the following recursive relation, obtained
by conditioning on the variables $X (v)$ for $|v| = 1$:
\begin{equation} \label{recurse}
\Psi (b_1 , \ldots , b_n ; D) = \left [ \int \Psi ({b_2 \over b_1} ,
   \cdots , {b_n \over b_1} ; D/v_1) \, d\mu (x_1) \right ]^{b_1}
\end{equation}
where for $D \subseteq \R^n$ and $(x_1 , \ldots , x_k) \in \R^k$,
the notation $D/(x_1 , \ldots , x_k)$ is used for the cross-section
of $D$ given by 
$$\{ (x_{k+1} , \ldots , x_n) \in \R^{n-k} \,:\, (x_1 , \ldots , 
   x_n) \in D \} .$$
Observe that the relation~(\ref{recurse}) together with the initial
condition $\Psi (b;D) = \mu(D)^b$ for $D \subseteq \R$ uniquely
determines $\Psi$ and in fact remains a valid inductive definition
when the arguments $b_i$ are positive reals, not necessarily integral.
Call this extension $\Psi$ as well, since it agrees with the old
$\Psi$ on integral arguments.  
We verify by induction on $n$ that for any $D \subseteq \R_n$, the 
function $\Psi (b_1 , \ldots , b_n ; D)$ is nonincreasing in each of
its $n$ arguments.  

This is clear for $n=1$ so fix $n > 1$ and observe the fact that
$\Psi$ is nonincreasing in $b_2 , \ldots , b_n$ follows directly
from the induction hypothesis and~(\ref{recurse}).  It remains
to check that $\Psi (b_1 , \ldots , b_n) \leq \Psi (b_1' , b_2 , 
\ldots , b_n)$ when $b_1 \geq b_1'$.  Rewriting this as
$$ \left [ \int \Psi ({b_2 \over b_1} ,
   \cdots , {b_n \over b_1} ; D/x_1) \, d\mu (x_1) \right ]^{b_1}
\leq \left [ \int \Psi ({b_2 \over b_1'} ,
   \cdots , {b_n \over b_1'} ; D/x_1) \, d\mu (x_1) \right ]^{b_1'}$$
we see it is just H\"older's inequality $\int h(x_1) \,d\mu (x_1)
\leq [\int h(x_1)^r \,d\mu (x_1)]^{1/r}$ with $r=b_1/b_1'$ 
applied to the function
$$h(x_1) = \left [ \int \Psi ({b_3 \over b_2} , \ldots , {b_n \over b_2};
   D/(x_1 , x_2)) \,d\mu (x_2) \right ]^{b_2 / b_1} .$$
(For $n=2$ take $h(x_1) = \mu (D/x_1)^{b_2/b_1}$).  

This proves Theorem~\ref{comparison} for spherically symmetric trees.
To obtain the general case requires the following lemma.

\begin{lem} \label{log convex}
For any Borel set $D \subseteq \R^n$, the function $\Psi (b_1 , 
\ldots , b_n ; D)$ defined by~(\ref{recurse}) is log-convex on
the positive orthant of $\R^n$.  
\end{lem}

\noindent{{\bf Proof}:}  Assume that $\mu^n (D) > 0$, since otherwise
$\Psi (\cdots ; D)$ is identically zero.  For $n=1$, $\ln \Psi
(\,\cdot\, ; D)$ is linear.  Proceeding by induction, fix $n>1$.
The relation~(\ref{recurse}) shows that $\ln \Psi (\cdots ; D)$
is homogeneous of degree 1:
$$\ln \Psi (\lambda b_1 , \ldots , \lambda b_n ; D) = 
  \lambda \ln \Psi (b_1 , \ldots , b_n ; D) .$$
Hence to prove convexity it suffices to verify that 
$\ln \Psi (1, b_2 , b_3 , \ldots , b_n ; D)$ is a convex
function of the positive variables $b_2 , \ldots , b_n$.  But
$$\Psi (1 , b_2  , \ldots , b_n ; D) = \int \Psi (b_2 , \ldots , b_n 
  ;D/x_1) \,d\mu (x_1)$$
and the sum or integral of log-convex functions is log-convex (see
for example \cite[p. 7-10]{Ar}).   $\Cox$

\noindent{{\bf Proof of Theorem~{\protect \ref{comparison}} completed}}: 
We verify
by induction that for any tree $T$ of depth $n$ with
generation cardinalities $|T_i| = b_i$ for $i=1 , \ldots , n$ and
for any $D \subseteq \R^n$, the probability $\phi (T;D)$ that
all paths $\rho , w_1 , \ldots , w_n$ in $T$ satisfy 
$(X(w_1) , \ldots , X(w_n)) \in D$ is at least $\Psi (b_1 ,
\ldots , b_n ; D)$, where as usual, $X(v)$ are i.i.d. with
common law $\mu$ and $\Psi$ is defined by~(\ref{recurse}).  

For every $w$ in the first generation $T_1$ of $T$,
let $T(w)$ denote the subtree $\{ \rho \} \cup \{ u \in T :
u \geq w \}$.  Let $b_i (w)$ denote the cardinality of the $i^{th}$
generation of $T(w)$.  (In particular, $b_1 (w) = b_0 (w) = 1$.)
Fix $D \subseteq \R^n$.  By the induction hypothesis
\begin{eqnarray*}
\phi (T(w) ; D) & \geq & \int \Psi (b_2 (w) , \ldots , b_n (w) ; D/x_1)
   \,d\mu (x_1) \\[2ex]
& = & \Psi (1 , b_2 (w) , \ldots , b_n (w); D)  .
\end{eqnarray*} 
Therefore, since $\sum_{w \in T_1} b_i (w) = b_i$ for $i \geq 1$,
using the log-convexity established in the previous lemma gives
\begin{eqnarray*}
\phi (T ; D) & = & \prod_{w \in T_1} \phi (T(w) ; D) \\[2ex]
& \geq & \prod_{w \in T_1} \Psi (1 , b_2 (w) , \ldots , b_n (w) ;D) \\[2ex]
& \geq & \Psi (1 , {b_2 \over b_1} , \ldots , {b_n \over b_1} ;D)^{b_1} \\[2ex]
& = & \Psi (b_1 , \ldots , b_n ; D)  .
\end{eqnarray*}
$\Cox$

Lemma~\ref{log convex} yields a simple description of domination 
between trees of height 2.  

\begin{pr} \label{height 2}
Let $\Gamma$ and $\Gamma'$ be trees of height 2.  For each vertex 
in $\Gamma_1$, count its children and order the numbers so 
obtained in a decreasing sequence $n_1 \geq n_2 \geq \cdots \geq n_b$.
Similarly obtain a sequence $n_1' \geq n_2' \geq \cdots \geq n_b'$
from $\Gamma'$, appending zeros to one of the sequences if necessary 
so that both have the same length $b$.  Then $\Gamma$ dominates $\Gamma'$
if and only if 
\begin{equation} \label{young}
\sum_{i > k} n_i \geq \sum_{i > k} n_i' \mbox{ for every } k \geq 0 .
\end{equation}
\end{pr}

\noindent{{\bf Remark}}

Thinking of the numbers $n_i$ as a partition of $\sum
n_i = |\Gamma_2|$, this is the order gotten by combining the usual
majorization order of partitions (in the reverse direction) with
the inclusion order (Young's lattice).  

\noindent{\bf Proof:}  We start by showing that the condition~(\ref{young})
implies domination.  Assume without loss of generality that 
\begin{equation} \label{same size}
\sum_{i=1}^b n_i = \sum_{i=1}^b n_i'
\end{equation}
since otherwise we could increase $n_1'$, thereby obtaining a tree
dominating $\Gamma'$, while condition~(\ref{young}) would remain
unaffected.  

Conditions~(\ref{young}) and~(\ref{same size}) together
imply that the vector $( n_i': i = 1, \ldots , b)$ {\em majorizes}
the vector  $( n_i : i = 1, \ldots , b)$  in the sense that
the latter is a convex combination of permutations of the former
(see \cite[Theorem 47]{HLP}).  To show that $\Gamma$ dominates $\Gamma'$ it
suffices to verify that for any $D \subseteq \R^2$
\begin{equation} \label{prod ht 2}
\prod_{i=1}^b \Psi (1,n_i ; D) \leq \prod_{i=1}^b \Psi (1,n_i' ; D) .
\end{equation}
Fixing $D$, let $h(n_1 , \ldots , n_b)$ denote the left-hand side
of~(\ref{prod ht 2}).  Clearly $h$ is invariant under permutations
of its arguments.  By Lemma~\ref{log convex} it is a product of 
log-convex functions and hence log-convex.  Since $(n_1 , \ldots , n_b)$
is a weighted average of permutations of $(n_1' , \ldots , n_b')$,
the inequality~(\ref{prod ht 2}) follows.

For the converse, assume $\Gamma$ dominates $\Gamma'$.  Let $r \geq 1$.
For $0 < \ee < 1$, take 
$$D_\ee = \left ( [0,\ee^r] \times [0,1] \right ) \, \cup \,  
   \left ( [0,1] \times [0,\ee] \right ) .$$  
If $\{ X(v) : v \in \Gamma \}$ are independent and uniform on $[0,1]$,
then the probability that $(X(v_1) , X(v_2)) \in D$ for every
path $\rho, v_1, v_2$ lies between
$\prod_{i=1}^b \ee^{\min (n_i , r)}$ and $\prod_{i=1}^b (2\ee)^{\min 
(n_i , r)}$.  When $\ee$ is sufficiently small, the assumption that
$\Gamma$ dominates $\Gamma'$ forces 
$$\sum_{i=1}^b \min (n_i , r) \geq \sum_{i=1}^b \min (n_i' , r) .$$
Choosing $r = n_k'$ yields
$$rk + \sum_{i>k} n_i \geq \sum_{i=1}^b \min (n_i' , r) = rk + \sum_{i>k}
   n_i' $$
proving~(\ref{young}).   $\Cox$

We end this section with some remarks and questions about domination.

\noindent{1. } Already for trees of height 3 the domination order is
somewhat mysterious.  Consider the trees $\Gamma$ and $\Gamma'$
in figure 1, where $\Gamma'$ is obtained from $\Gamma$ by gluing
together the vertices in the first generation.

\begin{picture}(300,100)
\put(80,80){\circle*{3}}
\put(60,60){\circle*{3}}
\put(100,60){\circle*{3}}
\put(60,40){\circle*{3}}
\put(85,40){\circle*{3}}
\put(115,40){\circle*{3}}
\put(45,20){\circle*{3}}
\put(60,20){\circle*{3}}
\put(75,20){\circle*{3}}
\put(85,20){\circle*{3}}
\put(115,20){\circle*{3}}
\put(80,80){\line(-1,-1){20}}
\put(80,80){\line(1,-1){20}}
\put(60,60){\line(0,-1){20}}
\put(100,60){\line(-3,-4){15}}
\put(100,60){\line(3,-4){15}}
\put(60,40){\line(-3,-4){15}}
\put(60,40){\line(0,-1){20}}
\put(60,40){\line(3,-4){15}}
\put(85,40){\line(0,-1){20}}
\put(115,40){\line(0,-1){20}}
\put(77,3){$\Gamma$}
\put(230,80){\circle*{3}}
\put(230,60){\circle*{3}}
\put(210,40){\circle*{3}}
\put(230,40){\circle*{3}}
\put(250,40){\circle*{3}}
\put(195,20){\circle*{3}}
\put(210,20){\circle*{3}}
\put(225,20){\circle*{3}}
\put(230,20){\circle*{3}}
\put(250,20){\circle*{3}}
\put(230,80){\line(0,-1){20}}
\put(230,60){\line(1,-1){20}}
\put(230,60){\line(-1,-1){20}}
\put(230,60){\line(0,-1){20}}
\put(210,40){\line(-3,-4){15}}
\put(210,40){\line(0,-1){20}}
\put(210,40){\line(3,-4){15}}
\put(230,40){\line(0,-1){20}}
\put(250,40){\line(0,-1){20}}
\put(230,3){$\Gamma'$}
\put(290,0){figure 1}
\end{picture}

\noindent{Intuitively}, it seems that $\Gamma$ should dominate $\Gamma'$,
but this is not the case.  Let 
$$D = \left ( [0,1/2] \times [0,1] \times [0,2/3] \right ) \, \cup \, 
\left ( [1/2,1] \times [0,1/2] \times [0,1] \right ) $$
and let $X(v)$ be uniform on $[0,1]$.  The probability that all paths
in $\Gamma$ have \\ $(X(v_1) , X(v_2) , X(v_3)) \in D$ is ${1075 \over 7776}$,
while the corresponding probability for $\Gamma'$ is only ${998 \over
7776}$.

\noindent{2. } To verify that a tree $\Gamma$ dominates another
tree $\Gamma'$, it suffices to consider the case in which the
i.i.d. variables $X(v)$ are uniform on $[0,1]$, since other variables
can be written as functions of these.  Along the same lines, 
Theorem~\ref{comparison} may be seen to sharpen a result from
\cite{BP1}.  There it was shown that any tree of height $n$ with
$b_n$ vertices in the $n^{th}$ generation is dominated by the tree
consisting of $b_n$ disjoint paths of length $n$ from the root.
The notion of domination in \cite{BP1} is for tree-indexed Markov chains;
we omit the easy proof that this is an equivalent definition of
domination.  

\noindent{3. } A counterexample of the type given in figure 1 cannot
occur if the set $D$ is restricted to being an upwardly closed
subset of $[0,1]^n$ (i.e. ${\bf x} \in D$ and ${\bf y} \geq {\bf x}$
coordinatewise imply ${\bf y} \in D$).  These subsets occur naturally
in percolation problems.  Domination for upwardly
closed sets subsumes gluing, as may be shown using the FKG inequality.  
Is this partial order any more tractable for trees of height greater than 2?

\section{Exponential transit times and nondecreasing \break growth functions}

\noindent{{\bf Proof of Theorem~{\protect \ref{increasing}}}}:
When $\sum_{n=1}^\infty f(n)^{-1} < \infty$, use the 
greedy algorithm to select the (a.s. unique) random path 
$\rho = v_0 , v_1 , v_2 , \ldots$ for which
$X(v(n))$ is minimal among $\{ X(w) : w \mbox{ is a child of }
v_{n-1} \}$.  Then 
$$\E X({v_n}) = \E [X({v_n}) \| v_{n-1}] = f(n)^{-1},$$
being the minimum of $f(n)$ standard exponentials.
Hence $\E \sum_{n=1}^\infty X({v_n}) < \infty$ so in particular,
$\sum_{n=1}^\infty X({v_n}) < \infty$ almost surely.

The more interesting case is when $\sum_{n=1}^\infty f(n)^{-1} = \infty$.
In this case, consider the weighted sums 
$$S^* (v) = \sum_{\rho < w \leq v} Y(w) $$
where $Y(w) = f(|w|) X(w)$.  Define 
$$M_n^* = \min_{|v| = n} S^* (v) .$$
For any vertex $w \neq \rho$ and any $\lambda > 0$ we have
$$\E \left [ e^{- \lambda Y(w)} \right ] = {1 \over 1 + \lambda f(|w|)} .$$
Multiplying these along a path yields 
$$\E \left [ e^{-\lambda S^* (v)} \right ] = \prod_{\rho < w \leq v}
   {1 \over 1 + \lambda f(|w|) } < \lambda^{-|v|} \prod_{j=1}^{|v|} 
   f(j)^{-1},$$
and summing over $\Gamma_n$ gives $\E \left [ \sum_{v \in \Gamma_n}
\exp (-\lambda S^* (v)) \right ] < \lambda^{-n}$.  Now for any
$\ee$ with $0 < \ee < 1$, Markov's inequality implies that
$$\P \left [ \exists v \in \Gamma_n : S^* (v) \leq {n \over \lambda}
  \ln ((1-\ee) \lambda) \right ] < (1-\ee)^n .$$
By Borel-Cantelli, it follows that $M_n^* \leq {n \over \lambda}
\ln ((1-\ee)\lambda)$ finitely often almost surely;
choosing $\lambda = e$ and letting $\ee \rightarrow 0$ leads
to 
\begin{equation} \label{deterministic}
\liminf_{n \rightarrow \infty} {1 \over n} M_n^* \geq e^{-1}
  \; \mbox{ almost surely.}  
\end{equation}

Deriving a lower estimate for $M_n$ from~(\ref{deterministic})
requires no probability theory.  Choose $\ee > 0$; with probability
one there exists an integer $N_\ee$ such that $M_k^* \geq k (e^{-1} - \ee)$
for all $k \geq N_\ee$.  For any path $\{ v_k : 0 \leq k \leq n \}$
starting from $v_0 = \rho$, summation by parts yields
\begin{eqnarray*}
\sum_{k=1}^n X(v_k) & = & \sum_{k=1}^n f(k)^{-1} [S^* (v_k) - S^*
   (v_{k-1})]  \\[2ex] 
& = & S^* (v_n) f(n+1)^{-1} + \sum_{k=1}^n S^* (v_k) [ f(k)^{-1}
   - f(k+1)^{-1} ] \\[2ex]
& > &  S^* (v_n) f(n+1)^{-1} + (e^{-1} - \ee) \sum_{k=1}^n k 
   [ f(k)^{-1} - f(k+1)^{-1} ] - C(N_\ee)
\end{eqnarray*}
where $C(N_\ee)$ depends only on $N_\ee$ and $f$.
The last inequality is the only place we use the assumption that
$f$ is nondecreasing.  Summing by parts again,
\begin{eqnarray*}
\sum_{k=1}^n X(v_k) & > & (e^{-1} - \ee) \sum_{j=1}^n f(j)^{-1} - 
   C(N_\ee) + f(n+1)^{-1} [ S^* (v_n ) - n (e^{-1} - \ee)] \\[2ex]
& \geq & (e^{-1} - \ee) \sum_{j=1}^n f(j)^{-1} -  C(N_\ee) .
\end{eqnarray*}
Thus $M_n$ is also greater than the right-hand side of the last 
inequality, which easily implies that almost surely
$$\liminf_{n \rightarrow \infty} M_n / \sum_{j=1}^n f(j)^{-1} 
   \geq e^{-1},$$
as long as the series $\sum f(j)^{-1}$ diverges.  This completes the proof
of the explosion criterion.  

To show that 
$$\limsup_{n \rightarrow \infty} M_n / \sum_{j=1}^n f(j)^{-1} 
   \leq e^{-1}$$
when the denominator tends to infinity and $f$ is unbounded
it suffices to exhibit, for every $\ee > 0$, an infinite path in
$\Gamma$ along which the condition
\begin{equation} \label{fail fo}
S(v) / \sum_{j=1}^{|v|} f(j)^{-1} \leq e^{-1} + \ee
\end{equation}
fails only finitely often.  This may be accomplished by a branching
process argument as in \cite{Pe}; the reader is referred there for 
greater detail.  

Let $k$ be any positive integer.  We prune the 
tree $\Gamma$ to obtain a random subtree $\Gamma'$ as follows.
First, take $\rho \in \Gamma'$; next, for $v \in \Gamma'$ and 
$w$ a child of $v$ in $\Gamma$, let $w \in \Gamma'$ if and only if
$X(w)$ has one of the $k$ least values among $\{ X(u) : u
\mbox{ is a child of $v$ in } \Gamma \}$.  Then $\Gamma'$ is a tree
which is eventually $k$-ary.  Again, let $Y(v) = f(|v|) X(v)$
with partial sums $S^*(v) = \sum_{\rho < w \leq v} Y(w)$ and
$M_n^*$ as before.  For $v \in \Gamma'$ the variables 
$\{ Y(w) : w \in \Gamma' \mbox{ is a child of } v \}$ are the
first $k$ order statistics of $f(|v|)$ independent exponentials
of mean $f(|v|)$.  As $|v| \rightarrow \infty$, we have assumed that
also $f(|v|) \rightarrow \infty$, so these order statistics are
converging weakly to the first $k$ hits of a mean 1 Poisson process.
[This is an easy consequence of the usual Poisson convergence theorem
and the fact that the common distribution of the original random variables
$X(v)$ has a density of one at the origin.]  In particular, for every
$\ee > 0$ and $v \in \Gamma'$, the joint distribution of 
$$\{ Y(w) : w \in \Gamma' , w \mbox{ is a child of } v \}$$
is stochastically dominated by  the first $k$ hits of a mean $1-\ee$ Poisson 
process.  In other words, conditional on the selection of the vertices
belonging to the subtree $\Gamma'$, we can couple the process $\{ Y(v) :
v \in \Gamma' \}$ to a process $\{ Z(v) : v \in \Gamma' \}$
for which the values of $Z$ over children of $v$ are independent 
for different $v$ and each distributed as a uniform random
permutation of the first $k$ order statistics of a mean $1-\ee$ Poisson 
process whenever $f(|v| +1) \geq k$; the coupling will satisfy 
$Z(v) > Y(v)$ for all $v$ such that $|v| > R$ for some $R$ depending
on $k$ and $\ee$.  

If $Z$ is a random variable distributed as a uniform random selection
of one of the first $k$ hits of a mean $1-\ee$ Poisson process, then 
$$\E e^{-\lambda Z} = k^{-1} \sum_{i=1}^k \left ( 1 + {\lambda \over
    1-\ee } \right )^{-i} = {1-\ee \over k \lambda} \left [ 1 - \left ( 1 +
    {\lambda \over 1-\ee} \right )^{-k} \right ] .$$
Choose $k = k(\ee)$ so that this is at least $(1-2\ee)/k\lambda$
for all $\lambda > 1/10$, say.  The rate function for $Z$ satisfies
\begin{eqnarray*}
m_Z (a) & \deq & \inf_{\lambda > 0} e^{\lambda a} \E e^{-\lambda Z} \\[2ex]
& \geq & {1 - 2\ee \over k} \inf_{\lambda > 0} e^{\lambda a}
   \min \{ {1 \over \lambda} , 10 \} \\[2ex]
& \geq & { 1-2\ee \over k} ea
\end{eqnarray*}
for $0 < a < 10/e$.  In particular if we fix $a=(1-3\ee)^{-1} e^{-1}$
then $m_Z (a) > k^{-1}$.  

Now the proof of (2.8) in \cite[p. 1237]{Pe} or the main 
result of \cite{Bi} shows that there is some integer $L$ for
which the following branching process is supercritical: $v \in \Gamma'$
begets $w \in \Gamma'$ if $|w| = |v| + L$ and $L^{-1} \sum_{v < u \leq w}
Z(u) \leq a$.  This implies that for such an $L$ there exists almost
surely a path $\rho , v_1 , v_2 , \ldots$ for which
\begin{equation} \label{stays close}
{1 \over L} \sum_{i = jL+1}^{jL+L} Y(v_i) \leq a \mbox{ for all but
 finitely many } j .
\end{equation}
From~(\ref{stays close}) we infer that 
$$\sum_{i=jL+1}^{jL+L} X(v_i) \leq aL\, f(jL+1)^{-1} \leq a
   \sum_{i=jL-L+1}^{jL}f(i)^{-1} $$
for all but finitely many $j$, and hence 
$$S(v_n) \leq a \sum_{i=1}^n f(i)^{-1} + O(1) .$$
Since $a = (1-3\ee)^{-1} e^{-1}$, letting $\ee \rightarrow 0$
finishes the proof.    $\Cox$
   
\section{Exponential transit times, arbitrary growth functions}

\noindent{{\bf Proof of Theorem~{\protect \ref{general}}}}:  
Consider first the case
$\sum \ft(n)^{-1} = \infty$.  If $\liminf_{n \rightarrow \infty}
( \prod_{i=1}^n f(i))^{1/n} < \infty$ then clearly there is no 
explosion since utilizing the gamma distribution shows that
$\P (M_n < C) \leq {C^n \over n!} \prod_{i=1}^n f(i)$.  Henceforth
we assume that 
$$\liminf_{n \rightarrow \infty} ( \prod_{i=1}^n f(i))^{1/n} = \infty ,$$ 
which is equivalent to $\ft (n) \rightarrow \infty$.  Under this
assumption, for every $n \geq 1$ there exists an $m \geq 1$ such that
$$\ft (n+1)^m \prod_{i=1}^n \ft (i) = \prod_{i=1}^{n+m} f(i) .$$
Consequently, there is an increasing sequence $\{ n(k) : k = 
1, 2, 3, \ldots \}$ for which 
$$\prod_{i=1}^{n(k)} \ft (i) = \prod_{i=1}^{n(k)} f (i) $$
for all $k$.  

To establish the lower bound on $\limsup M_n / \sum_{i=1}^n \ft 
(i)^{-1}$, and hence also that no explosion occurs, use the
proof of Theorem~\ref{increasing} but with $Y(v) = \ft (|v|) X(v)$.
As in the proof of that theorem, the partial sums $S^* (v) 
= \sum_{\rho < w \leq v} Y(w)$ satisfy 
$$\E \left [ e^{-\lambda S^* (v)} \right ] \leq \lambda^{-|v|}
   \prod_{i=1}^{|v|} \ft (i) $$
for $\lambda > 0$.  When $|v| = n(k)$, we can substitute $f$
for $\ft$ on the right-hand side.  It follows for $\ee > 0$,
that with probability one the inequality 
$$M_{n(k)}^* \leq {n(k) \over \lambda} \ln ((1-\ee)\lambda) $$
holds for only finitely many $k$.  Taking $\lambda = e$ we get,
as before,
$$\liminf_{k \rightarrow \infty} M_{n(k)} / \sum_{i=1}^{n(k)}
  \ft(i)^{-1} \geq e^{-1}$$
so in particular
$$\limsup_{n \rightarrow \infty} M_n / \sum_{i=1}^n
  \ft(i)^{-1} \geq e^{-1} .$$

Next, let us bound this $\limsup$ from above, still assuming that
$\sum \ft (n)^{-1} = \infty$ and $\ft (n) \rightarrow \infty$.
Let $g(n)$ denote the integer part of $\ft(n)$ and let
$\Gamma'$ be a spherically symmetric tree with growth function
$g$.  And $\iid$ random variables $\{ X (w) \}$.
Apply Theorem~\ref{increasing} to $\Gamma'$.  Examining the
proof at~(\ref{fail fo}), we see that for any $\ee > 0$ there 
exists almost surely some path $\rho , v_1 , v_2 , \ldots$ satisfying
\begin{equation} \label{approximant}
X (v_n) \leq (e^{-1} + \ee) \sum_{i=1}^n g(i)^{-1}
\end{equation}
for all but finitely many $n$.  From
Theorem~\ref{comparison} it follows that a path satisfying~(\ref{approximant}) 
exists almost surely in $\Gamma$ as well as $\Gamma'$.
Since $\sum_{i=1}^n g(i)^{-1} / \sum_{i=1}^n \ft (i)^{-1}
\rightarrow 1$, this implies that almost surely 
$$\limsup_{n \rightarrow \infty} M_n / \sum_{i=1}^n
  \ft(i)^{-1} \leq e^{-1}. $$

Finally, in the case when $\sum \ft (n)^{-1}$ converges, define
$g$ and $\Gamma'$ as above and apply Theorem~\ref{increasing}
to conclude that for some fixed $L > 0$
$$ \liminf_{n \rightarrow \infty} \P \left [ \exists w \in \Gamma_n' 
   : S'(w) \leq L \right ] > 0 .$$
This together with Theorem~\ref{comparison} and the zero-one law 
for explosions proves that the probability of an explosion is 1.
$\Cox$

\section{Other transit-time distributions}

Theorems~\ref{increasing} and~\ref{general} hold without change
when the i.i.d. transit times $X(v)$ have any distribution function
$G$ for which $G(t)/t = 1+o(1)$ as $t \rightarrow 0$.  More generally
we have 

\begin{pr} \label{power law}
Let $\Gamma$ be an infinite spherically symmetric tree with growth
function $f$.  Suppose $\{ X(v) : v \in \Gamma \}$ are i.i.d. with
their common distribution function $G$ satisfying
$$\lim_{t \rightarrow 0} G(t) t^{-\alpha} = c > 0 $$
for some $\alpha > 0$.
With $S(v)$ and $M_n$ defined as throughout, we then have the
following analogues of Theorems~\ref{increasing} and~\ref{general}.  
\begin{quotation}
$(i)$ If $f$ is nondecreasing then there is an explosion
($\sup M_n < \infty$) if and only if $\sum f(n)^{-1/\alpha}$
converges.  If the sum diverges then
$$\lim_{n \rightarrow \infty} M_n / \sum_{j=1}^n f(j)^{-1/\alpha}
  = \alpha e^{-1} [ c \Gamma (1+\alpha)]^{-1/\alpha} .$$ 

$(ii)$ For a general growth function $f$, define $\ft$ as in 
the preface to Theorem~\ref{general}.  There is almost surely an 
explosion if and only if $\sum \ft (n)^{-1/\alpha}$ converges.
If the sum diverges then with probability one,
$$\limsup_{n \rightarrow \infty} M_n / \sum_{j=1}^n \ft (j)^{-1/\alpha}
  = \alpha e^{-1} [ c \Gamma (1+\alpha)]^{-1/\alpha} .$$ 
\end{quotation}
\end{pr}

\noindent{\bf Proof:}  We only discuss the modifications needed in the
proof for the exponential case.  For $(i)$, first assume that 
$\sum f(n)^{-1/\alpha}$ converges.  For every vertex $v \in \Gamma_{n-1}$,
the number of its children $w \in \Gamma_n$ for which 
\begin{equation} \label{six}
X(w) \leq \left [ {c \over 2} f(n) \right ]^{-1/\alpha}
\end{equation}
converges in distribution to a Poisson with mean 2.  Comparing to 
a Galton-Watson process shows that~(\ref{six}) holds along an
entire infinite path with positive probability and this implies
almost sure explosion.  

Next, assume that $\sum f(n)^{-1/\alpha} = \infty$ and that
$f(n) \rightarrow \infty$.  Given $\ee > 0$, choose $\dd > 0$ so that
$G(t) < (c+\ee)t^\alpha$ for $0 < t < \dd$.  Estimate the moment 
generating function 
\begin{eqnarray*} 
\E [e^{-\lambda X(v)} ] & = & \int_0^\infty e^{-\lambda t} \,dG(t) \\[2ex]
& = & \int_0^\infty G(t) \lambda e^{-\lambda t} \,dt \\[2ex]
& \leq & \int_0^\dd (c+\ee) t^\alpha \lambda e^{-\lambda t} \, dt
   + e^{-\lambda \dd} \\[2ex]
& \leq & (c+\ee) \lambda^{-\alpha} \Gamma (1+\alpha) + 
   e^{-\lambda \dd} \\[2ex]
& \leq & (c+2\ee) \lambda^{-\alpha} \Gamma (1+\alpha)
\end{eqnarray*}
for large positive $\lambda$.  Define $S^* (v) = \sum_{\rho < w \leq v} 
f(|w|)^{1/\alpha} X(w)$ and $M_n^* = \min_{|v| = n} S^* (v)$.
Since $f(n) \rightarrow \infty$ it follows that
$$\E [ \exp (-\lambda S^* (v))] \leq A(\lambda) \prod_{j=1}^{|v|}
   (c+2\ee) \lambda^{-\alpha} f(j)^{-1} \Gamma (1+\alpha) $$
where $A(\lambda) > 0$, and hence 
$$\E \left [ \sum_{|v|=n} \exp (-\lambda S^* (v)) \right ] \leq A(\lambda) 
  \left [ (c+2\ee) \lambda^{-\alpha} \Gamma (1+\alpha) \right ]^n . $$
As in the proof of Theorem~\ref{increasing}, optimizing over $\lambda$
yields 
$$\liminf {1 \over n} M_n^* \geq e^{-1} \alpha [c \Gamma (1+\alpha
   )]^{-1/\alpha}$$
and consequently, with probability one,
$$\liminf_{n \rightarrow \infty} M_n / \sum_{j=1}^n f(j)^{-1/\alpha}
  \geq e^{-1} [c \Gamma (1+\alpha)]^{-1/\alpha} .$$
The rest of the proof proceeds as in Theorem~\ref{increasing}.

Finally, the only change in the proof of $(ii)$ is to note that the
operation $f \mapsto \ft$ commutes with taking powers.  

\noindent{\bf Question:} is there a simple explosion criterion for 
arbitrary transit time distributions?  At least when the
growth function $f$ is nondecreasing and $G$ is strictly monotone
and continuous, it seems possible that an explosion occurs if and only if
\begin{equation} \label{criterion}
\sum_{n=1}^\infty G^{-1} \left ( { 1 \over f(n)} \right ) < \infty 
\end{equation}
where $G^{-1}$ is the inverse function to $G$.
The technique used to obtain Proposition~\ref{power law} is powerful
enough to verify this criterion for a slightly more general class
of distributions, but the general case has eluded us.

{\bf ACKNOWLEDGEMENTS}

\noindent{We} are grateful to Enrique Andjel for posing the problem 
which led to this work.  The second author would also like to thank
Itai Benjamini for useful discussions and J.F. Mela and F. Parreau 
for their hospitality at the Universit\'e Paris-Nord.

\renewcommand{\baselinestretch}{1.0}\large

\end{document}